\PassOptionsToPackage{hidelinks}{hyperref}
\documentclass[submission,nomarks,a4paper]{dmtcs-episciences}
\usepackage[T1]{fontenc}
\usepackage[utf8]{inputenc}
\usepackage{amsmath,amssymb,mathtools}
\usepackage{booktabs,array,tabularx,microtype,xurl,placeins}
\usepackage[numbers,sort&compress]{natbib}
\numberwithin{equation}{section}
\newtheorem{theorem}{Theorem}[section]
\newtheorem{proposition}[theorem]{Proposition}
\newtheorem{lemma}[theorem]{Lemma}
\newtheorem{remark}[theorem]{Remark}
\hypersetup{pdftitle={Fixed-level profiles and fluctuations in age-dependent recursive trees},
 pdfauthor={Marek Gałązka and Jerzy Szymański}}

\newcommand{\E}{\mathbb E}
\newcommand{\Pp}{\mathbb P}
\newcommand{\F}{\mathcal F}
\newcommand{\1}{\mathbf 1}
\newcommand{\Var}{\operatorname{Var}}
\newcommand{\Cov}{\operatorname{Cov}}
\newcommand{\norm}[1]{\left\lVert#1\right\rVert_2}
\newcommand{\dto}{\xrightarrow{\mathrm d}}
\newcommand{\Ltwo}{\xrightarrow{L^2}}
\newcommand{\doi}[1]{\href{https://doi.org/#1}{doi:\nolinkurl{#1}}}

\title[Fixed-level profiles in age-dependent recursive trees]{Fixed-level profiles and fluctuations in age-dependent recursive trees}
\author[Marek Ga\l\k{a}zka and Jerzy Szyma\'nski]{Marek Ga\l\k{a}zka\thanks{Corresponding author: \texttt{galazka@amu.edu.pl}.} \and Jerzy Szyma\'nski\thanks{\texttt{jesz@amu.edu.pl}.}}
\affiliation{Department of Algorithm Theory and Data Security, Faculty of Mathematics and Computer Science, Adam Mickiewicz University, Pozna\'n, Poland}
\keywords{random recursive trees, age-dependent attachment, level profile, covariance kernel, stabilization, multivariate central limit theorem}
\makeatletter
\AtBeginDocument{\gdef\@firsthead{}\gdef\@copyright{}}
\makeatother
\begin{document}
\maketitle
\begin{abstract}
How many vertices of a growing recursive tree lie at a prescribed distance from the root? We compare two linear attachment rules. When an arriving vertex chooses a parent with probability proportional to its label, each fixed level eventually stops growing. We determine the limiting covariance kernel and the probability that a level will change again after time $n$. Writing $L_\infty(h)$ for its final size, we show that $2^{-h}L_\infty(h)$ converges in mean square as $h$ tends to infinity to a random variable with mean $1$ and variance $1/3$. Consequently the limiting variance of the level size is asymptotic to $4^h/3$. When parents are chosen with probability proportional to their age, each fixed level grows. Its centered fluctuations have a joint Gaussian limit whose covariance constants come from repeated logarithmic integration of attachment noise at the root. An exact finite-size calculation illustrates the difference between the limiting and finite-size variances.
\end{abstract}
\noindent\textbf{Mathematics Subject Classification (2020).} Primary 60C05. Secondary 05C80, 60F05.

\section{Introduction}\label{sec:intro}
The first few levels of a recursive tree describe the neighborhood of its oldest vertex. Their sizes depend on how often new vertices attach to an early ancestor. Two trees can have the same insertion-depth distribution, and hence the same expected profile, without having the same fluctuations across levels. Here the question is how a linear preference for the age of a parent changes these fluctuations.

Vertices arrive in order $1,2,\ldots$. At time $n$, model $N$ assigns weight $i$ to a possible parent $i<n$, and model $O$ assigns weight $n-i$. In model $N$ the probabilities of attachment to the root sum to two over the life of the tree. In model $O$ their partial sum grows as $2\log n$. Both rules belong to the nonuniform recursive-tree setting of Szyma\'nski~\cite{Szymanski1987}. Hofri and Mahmoud~\cite[Section 10.2]{HofriMahmoud} analyze age affinity in model $O$, including individual degrees and insertion depth. Ga\l\k{a}zka et al.~\cite{Galazka2020} study degree characteristics under both label-based rules. For model $N$, Lyon and Mahmoud~\cite{Lyon2020} identify the insertion-depth law and the variance of the root degree. Their depth law also occurs for binary search trees. Model $N$ is the exponent-one member of their power-weight family~\cite{Lyon2022}.

We write $L_n(h)$ for the number of vertices at depth $h$, $W_n(h)$ for the sum of their labels and $S_n(h)$ for the sum of their undirected degrees. Table~\ref{tab:scope} records what can be recovered from the earlier marginal results and where a joint argument enters. For example, $S_n(0)=L_n(1)$, whereas the variance of $S_n(1)$ also involves the covariance of two different level counts.

\begin{table}[!htbp]
\centering\small
\caption{Earlier inputs and results established here. A marginal depth law determines both $\E L_n(h)$ and $\E W_n(h)$ by summation.}\label{tab:scope}
\begin{tabularx}{\linewidth}{@{}l>{\raggedright\arraybackslash}X>{\raggedright\arraybackslash}X@{}}
\toprule
 & Earlier inputs or direct consequences & Role in this paper\\
\midrule
$N$, $L$ & Insertion-depth law, fixed-level means and eventual stabilization \cite{Lyon2020}. & Joint $L^2$ limits, kernel $K$, its growth with depth and sharp late-increment probabilities.\\
$N$, $W$ & Its mean is a label-weighted sum of the same marginal probabilities. & Conditional predictor and mixed second moments that determine $K$.\\
$N$, $S$ & Finite root variance and its limit $K(1,1)$ \cite[Prop.~4.1]{Lyon2020}. & Joint degree-sum covariances from adjacent levels, including the variance at depth one.\\
$O$, $L$ & Age-affinity rule and marginal depth analysis \cite[Sec.~10.2]{HofriMahmoud}. & Fixed-level covariance constants and a joint Gaussian limit after exact centering.\\
$O$, $W$ & Conditional attachment weights are linear in the parent label. & Moment estimates control the label-weighted term in the centered recurrence.\\
$O$, $S$ & Individual-vertex degree results \cite{HofriMahmoud,Galazka2020} concern a different statistic. & Joint Gaussian limit and covariance constants for sums over whole levels.\\
\bottomrule
\end{tabularx}
\end{table}

\paragraph{Related models.}
For uniform recursive trees, Drmota and Hwang~\cite{DrmotaHwang2005} and Fuchs, Hwang and Neininger~\cite{Fuchs2006} analyze profile fluctuations. Iksanov and Kabluchko~\cite{Iksanov2018} prove a functional limit for fixed levels using a continuous-time branching process. Kabluchko, Marynych and Sulzbach~\cite{KMS2017} derive profile expansions in a branching-random-walk framework that includes the uniform tree and its logarithmic bulk. For weighted recursive trees, S\'enizergues~\cite{Senizergues2021} studies the profile and geometry for prescribed vertex weights, Pain and S\'enizergues~\cite{Pain2022} study height, and Lodewijks and Ortgiese~\cite{Lodewijks2024} study degree extremes for random vertex weights. Model $N$ has prescribed weights $w_i=i$. In model $O$ the odds for two fixed potential parents vary with the insertion time. Section~\ref{sec:clock} gives the resulting obstruction to applying these theorems directly.

Scaled attachment gives another comparison. The parent label divided by $n$ has limiting density $2x$ in model $N$ and $2(1-x)$ in model $O$. Devroye, Fawzi and Fraiman~\cite{Devroye2010,Devroye2012} use scaled attachment to analyze typical and extremal depths. Exponential label weights lead to a different family studied by Aguech et al.~\cite{Aguech2025}. Kuba and Panholzer~\cite{Kuba2007} study the degree of a vertex with a specified label in increasing trees. Here degrees are summed over vertices selected by depth.

\paragraph{Contributions.}
For model $N$ we obtain a convergent series for the mixed limiting covariances and sharp probabilities for late changes. In a second limit, taken after the tree has grown to infinity, the same martingale gives the covariance asymptotic $K(h,k)\sim 2^{h+k}/3$ as both depths tend to infinity. For model $O$ we isolate the label-weighted correction and prove that root attachments generate the joint fluctuations. This yields the uniform tree's limiting correlations with explicit, different covariance constants. Section~\ref{sec:main} states the results, Sections~\ref{sec:newer} and~\ref{sec:older} give the proofs, and Section~\ref{sec:simulation} compares their predictions with exact finite-size moments.

\section{Models and the uniform benchmark}\label{sec:models}
\subsection{Attachment rules and notation}
The tree begins with the root labeled $1$. For $n\ge2$, vertex $n$ chooses a parent $\Gamma_n\in\{1,\ldots,n-1\}$. The random variables $\Gamma_2,\Gamma_3,\ldots$ are independent. The two parent distributions are
\begin{align}
 p^{N}_{ni}&=\Pp(\Gamma_n=i)=\frac{2i}{n(n-1)},\label{eq:pN}\\
 p^{O}_{ni}&=\Pp(\Gamma_n=i)=\frac{2(n-i)}{n(n-1)}.\label{eq:pO}
\end{align}
The letters $N$ and $O$ denote newer- and older-vertex preference. Age here means insertion rank: at time $n$, vertex $i$ has age $n-i$. Let $D(v)$ be the distance of vertex $v$ from the root, and let $\F_n=\sigma(\Gamma_2,\ldots,\Gamma_n)$. For $h\ge0$, put
\begin{align}
 I_j(h)&=\1_{\{D(j)=h\}},\label{eq:indicator}\\
 L_n(h)&=\sum_{j=1}^n I_j(h),\label{eq:profile}\\
 W_n(h)&=\sum_{j=1}^n jI_j(h),\label{eq:weight}\\
 S_n(h)&=\sum_{v:D(v)=h}\deg_n(v).\label{eq:degree}
\end{align}
Degrees are undirected, and the root has degree zero when $n=1$. Thus $L_n(0)=W_n(0)=1$. A superscript specifies the model when needed. We also write
\begin{equation}\label{eq:marginal}
 \ell_n(h)=\Pp(D(n)=h),\qquad F_n(h)=\E W_n(h).
\end{equation}
For the bounds and finite-size formulas we use, throughout,
\[
 \tau_n=1+\log n,\qquad H_m=\sum_{j=1}^m\frac1j,\qquad
 H_m^{(q)}=\sum_{j=1}^m j^{-q}\quad(q\ge1).
\]
In particular $H_m^{(2)}$ is the second-order harmonic number. Empty sums and quantities indexed by a negative level are zero. Constants in fixed-level bounds may depend on the indicated levels.

\begin{lemma}[Profile identities]\label{lem:identities}
For either model, $n\ge1$ and $h\ge0$,
\begin{align}
 S_n(h)&=L_n(h)+L_n(h+1)-\1_{\{h=0\}},\label{eq:identity}\\
 \E L_n(h)&=\sum_{j=1}^n\ell_j(h).\label{eq:summeans}
\end{align}
For $n\ge2$ and $h\ge1$,
\begin{equation}\label{eq:ellrec}
 \ell_n(h)=\sum_{i<n}p_{ni}\ell_i(h-1).
\end{equation}
The conditional depth probabilities are
\begin{align}
 \pi_n^N(h)&:=\Pp(D(n)=h\mid\F_{n-1})
       =\frac{2W_{n-1}^N(h-1)}{n(n-1)},\label{eq:piN}\\
 \pi_n^O(h)&:=\Pp(D(n)=h\mid\F_{n-1})
       =\frac{2L_{n-1}^O(h-1)}{n-1}-\frac{2W_{n-1}^O(h-1)}{n(n-1)}.\label{eq:piO}
\end{align}
\end{lemma}
\begin{proof}
Each non-root vertex at level $h$ has one parent edge. The children of level $h$ are exactly the vertices at level $h+1$, so counting incident edges gives~\eqref{eq:identity}. Summing the indicators gives~\eqref{eq:summeans}. The depth of a possible parent $i<n$ is measurable with respect to $\F_{n-1}$. Independence of $\Gamma_n$ from that sigma-field therefore gives~\eqref{eq:ellrec} and, on summing the attachment probabilities over the appropriate level,~\eqref{eq:piN}--\eqref{eq:piO}.
\end{proof}

\subsection{Classical uniform recursive trees}\label{sec:uniform}
For comparison, let $U$ denote the model with $p^U_{ni}=1/(n-1)$. For a fixed positive integer $r$, its exact mean is
\begin{equation}\label{eq:umean}
 \E L_n^U(r)=[z^r]\prod_{j=1}^{n-1}\left(1+\frac zj\right)
 =\frac{(\log n)^r}{r!}+O(\tau_n^{r-1}).
\end{equation}
Here $[z^r]$ means coefficient extraction. The product follows by taking expectations in the one-step profile recurrence. The error estimate follows by expanding $H_{n-1}^r$ and bounding terms with a repeated index.

The fixed-level results in~\cite{DrmotaHwang2005,Fuchs2006,Iksanov2018}, in our normalization, give
\begin{equation}\label{eq:ucov}
 \Cov(L_n^U(r),L_n^U(s))\sim
 \frac{(\log n)^{r+s-1}}{(r-1)!(s-1)!(r+s-1)},\qquad r,s\ge1.
\end{equation}
For every fixed $m\ge1$,
\begin{equation}\label{eq:uclt}
 \left(\frac{L_n^U(r)-\E L_n^U(r)}{(\log n)^{r-1/2}}\right)_{r=1}^m
 \dto (G_r^U)_{r=1}^m,
\end{equation}
where, for a standard Brownian motion $B$,
\begin{equation}\label{eq:ubrown}
 G_r^U=\frac{1}{(r-1)!}\int_0^1(1-t)^{r-1}\,dB(t).
\end{equation}
The covariance of this Gaussian vector is the constant in~\eqref{eq:ucov}. The mean estimate~\eqref{eq:umean} justifies replacing leading-order centering in the cited limit theorem by exact expectations. Appendix~\ref{app:uniform} also gives a direct second-moment derivation of~\eqref{eq:ucov}, so the covariance asymptotic is not inferred from weak convergence alone.

By~\eqref{eq:identity}, for fixed $h,k\ge0$,
\begin{align}
 \E S_n^U(h)&\sim\frac{(\log n)^{h+1}}{(h+1)!},\label{eq:usmean}\\
 \Cov(S_n^U(h),S_n^U(k))&\sim
 \frac{(\log n)^{h+k+1}}{h!k!(h+k+1)}.\label{eq:uscov}
\end{align}
The centered degree-sum vector, with coordinate $h$ divided by $(\log n)^{h+1/2}$, has Gaussian limit $(G_{h+1}^U)_h$. The same Brownian motion drives all the limiting coordinates.

\section{Main results}\label{sec:main}
\subsection{Newer-vertex preference}
For integers $r,s\ge0$, define the deterministic second moments
\begin{equation}\label{eq:Vdef}
 V_n(r,s)=\E[W_n^N(r)W_n^N(s)].
\end{equation}
They can be computed recursively by Lemma~\ref{lem:weighted}. Write $\zeta$ for the Riemann zeta function.

\begin{theorem}[Joint limits in model $N$]\label{thm:newer}
Every finite vector of the variables $L_n^N(h)$ and $S_n^N(h)$ converges in $L^2$ to a finite random vector. Each coordinate is eventually constant almost surely. Its means are
\begin{align}
 \E L_\infty^N(h)&=2^h,\label{eq:Nmeanlim}\\
 \E S_\infty^N(h)&=3\cdot2^h-\1_{\{h=0\}}.\label{eq:Nsmeanlim}
\end{align}
For $h,k\ge1$, put $K(h,k)=\Cov(L_\infty^N(h),L_\infty^N(k))$. Then
\begin{align}
 K(h,k)&=2^{h+k}\bigl(1-2^{-\min(h,k)}\bigr)\notag\\
 &\quad-4\sum_{m=1}^\infty\frac{1}{m^2(m+1)^2}
     \sum_{r=0}^{h-1}\sum_{s=0}^{k-1}2^{h+k-2-r-s}V_m(r,s).\label{eq:Kseries}
\end{align}
The series is absolutely convergent. Its tail after $m=M$ is
\begin{equation}\label{eq:Ktail}
 \begin{cases}
 O\bigl(\tau_M^{\min(h,k)-2}/M^2\bigr),&h,k\ge2,\\
 O\bigl(\tau_M^{\max(h,k)-1}/M^3\bigr),&\min(h,k)=1.
 \end{cases}
\end{equation}
Set $K(0,k)=K(h,0)=0$. The degree-sum covariances are
\begin{equation}\label{eq:Kdegree}
 \Cov(S_\infty^N(h),S_\infty^N(k))
 =\sum_{a=0}^1\sum_{b=0}^1K(h+a,k+b).
\end{equation}
In particular,
\begin{align}
 K(1,1)&=14-8\zeta(2),\label{eq:K11}\\
 K(1,2)&=60-24\zeta(2)-16\zeta(3),\label{eq:K12}\\
 K(2,2)&=340-80\zeta(2)-112\zeta(3)-64\zeta(4).\label{eq:K22}
\end{align}
Consequently,
\begin{equation}\label{eq:Nsvar}
 \Var(S_\infty^N(1))=474-136\zeta(2)-144\zeta(3)-64\zeta(4).
\end{equation}
\end{theorem}
The first moments and almost-sure stabilization are preparatory consequences of the known marginal law, as explained in Section~\ref{sec:Nfirst}. Formula~\eqref{eq:K11} is the root-variance limit in~\cite[Proposition 4.1]{Lyon2020}, numerically $0.8405274652$. The constants in~\eqref{eq:K12}, \eqref{eq:K22} and~\eqref{eq:Nsvar} are approximately $1.288672$, $4.506215$ and $7.924086$, respectively.

\begin{theorem}[Sharp stabilization errors]\label{thm:tails}
In model $N$ write $\Delta^L_n(h)=L_\infty^N(h)-L_n^N(h)$, $\Delta^S_n(h)=S_\infty^N(h)-S_n^N(h)$ and $b_n(d)=\E\Delta^L_n(d)$. Take $(\Delta,d)=(\Delta^L_n(h),h)$ for fixed $h\ge1$ or $(\Delta,d)=(\Delta^S_n(h),h+1)$ for fixed $h\ge0$. In both cases,
\begin{equation}\label{eq:jointtail}
 \E[\Delta^2]\sim\Pp(\Delta>0)\sim\E\Delta
 \sim b_n(d)\sim\frac{2^d(\log n)^{d-1}}{(d-1)!\,n}.
\end{equation}
Moreover $\Pp(\Delta=1\mid\Delta>0)\to1$. The root count has $\Delta^L_n(0)=0$, and $\Delta^S_n(0)=\Delta^L_n(1)$ with $b_n(1)=2/n$.
\end{theorem}

\subsection{Older-vertex preference}
\begin{theorem}[Means, covariances and Gaussian limits in model $O$]\label{thm:older}
For fixed $r,s\ge1$,
\begin{align}
 \E L_n^O(r)&\sim\frac{2^r}{r!}(\log n)^r,\label{eq:Omean}\\
 \Cov(L_n^O(r),L_n^O(s))&\sim C_{r,s}(\log n)^{r+s-1},\label{eq:Ocov}
\end{align}
where
\begin{equation}\label{eq:C}
 C_{r,s}=\frac{2^{r+s-1}}{(r-1)!(s-1)!(r+s-1)}.
\end{equation}
For every fixed $m\ge1$,
\begin{equation}\label{eq:Oclt}
 \left(\frac{L_n^O(r)-\E L_n^O(r)}{(\log n)^{r-1/2}}\right)_{r=1}^m
 \dto\mathcal N_m\bigl(0,(C_{r,s})_{r,s=1}^m\bigr).
\end{equation}
For fixed $h,k\ge0$, the corresponding degree-sum statements are
\begin{align}
 \E S_n^O(h)&\sim\frac{2^{h+1}}{(h+1)!}(\log n)^{h+1},\label{eq:Osmean}\\
 \Cov(S_n^O(h),S_n^O(k))&\sim D_{h,k}(\log n)^{h+k+1},\label{eq:Oscov}
\end{align}
with
\begin{equation}\label{eq:D}
 D_{h,k}=\frac{2^{h+k+1}}{h!k!(h+k+1)}.
\end{equation}
For every fixed $d\ge0$,
\begin{equation}\label{eq:Osclt}
 \left(\frac{S_n^O(h)-\E S_n^O(h)}{(\log n)^{h+1/2}}\right)_{h=0}^d
 \dto\mathcal N_{d+1}\bigl(0,(D_{h,k})_{h,k=0}^d\bigr).
\end{equation}
Both Gaussian covariance matrices are positive definite. In particular,
\begin{equation}\label{eq:Osvar}
 \Var(S_n^O(h))\sim\frac{2^{2h+1}}{(h!)^2(2h+1)}(\log n)^{2h+1}.
\end{equation}
The normalized uncentered quantities satisfy
\begin{align}
 \frac{L_n^O(r)}{(\log n)^r}&\Ltwo\frac{2^r}{r!},\label{eq:OL2}\\
 \frac{S_n^O(h)}{(\log n)^{h+1}}&\Ltwo\frac{2^{h+1}}{(h+1)!}.\label{eq:OSL2}
\end{align}
\end{theorem}
All Gaussian statements here use exact expectations for centering. The limits can be represented using one standard Brownian motion:
\begin{equation}\label{eq:Obrown}
 G_r^O=\frac{2^{r-1}\sqrt2}{(r-1)!}\int_0^1(1-t)^{r-1}\,dB(t).
\end{equation}
Thus the normalized limiting correlations agree with those of the uniform model, although the means and covariance constants differ. In particular, for $r,s\ge1$,
\begin{equation}\label{eq:corr}
 \lim_{n\to\infty}\operatorname{Corr}(L_n^O(r),L_n^O(s))
 =\frac{\sqrt{(2r-1)(2s-1)}}{r+s-1}.
\end{equation}

\subsection{The logarithmic clock in model \texorpdfstring{$O$}{O}}\label{sec:clock}
For a fixed parent $i$, the attachment probability at arrival $n$ is $(2/n)(1-(i-1)/(n-1))$, exactly $2/n$ at the root. In the uniform tree it is $1/(n-1)$. Thus the root noise accumulates variance $2\log n+O(1)$, and each subsequent leading drift integrates against $2\,d\log n$. There are $r-1$ such integrations for level $r$. In a covariance between levels $r$ and $s$ the two propagated weights and the common innovation variance give $2^{r-1}2^{s-1}2=2^{r+s-1}$. This explains the constant in~\eqref{eq:C} and the cancellation in~\eqref{eq:corr}.

A common change of clock leaves the relative rates of two competing parents unchanged. For fixed $i<j<n$ the odds in model $O$ are $(n-i)/(n-j)$, which vary with $n$. At the level of the profile, the exact conditional drift in~\eqref{eq:piO} depends on the sum of the parent labels $W_{n-1}(r-1)$ as well as on the number of parents. A bound for this second term requires an estimate beyond both the uniform branching-process embedding and the results on weighted recursive trees~\cite{Senizergues2021,Pain2022,Lodewijks2024} recalled in Section~\ref{sec:intro}. Lemmas~\ref{lem:moments} and~\ref{lem:rootapprox} show that after centering, dropping the label correction and higher-level insertion noise costs $O(\tau_n^{r-1})$ in $L^2$ for $r\ge2$. That estimate justifies the doubled logarithmic clock at the fluctuation scale.

\subsection{Growth of the limiting covariance kernel}\label{sec:large}
\begin{proposition}[Large depths in the infinite tree]\label{prop:large}
In model $N$ there is a random variable $Y$ with $\E Y=1$, $\kappa:=\Var Y=1/3$ and $\E Y^2=4/3$ such that
\begin{equation}\label{eq:deepL}
 \norm{2^{-h}L_\infty^N(h)-Y}\le 2^{-h/2},\qquad h\ge1.
\end{equation}
Consequently, as $\min(h,k)\to\infty$,
\begin{equation}\label{eq:deepK}
 \frac{K(h,k)}{2^{h+k}}\longrightarrow\frac13,\qquad
 \Var L_\infty^N(h)\sim\frac{4^h}{3}.
\end{equation}
\end{proposition}
Here the infinite-tree limit is taken at each fixed level first. This order of limits is distinct from following a level $h=h(n)$ in a finite tree. Section~\ref{sec:deepproof} proves the proposition and gives numerical enclosures for three further entries of $K$.

\begin{remark}[Exploring a fixed number of levels]\label{rem:work}
Expanding the adjacency list of each vertex through depth $h$ inspects $Q_n(h)=\sum_{r=0}^h S_n(r)=2\sum_{r=1}^hL_n(r)+L_n(h+1)$ entries. In model $N$ this cost stabilizes in $L^2$ and almost surely, with mean limit $3\cdot2^{h+1}-4$. In model $O$ only $L_n(h+1)$ contributes at the leading scales. Thus $Q_n(h)/(\log n)^{h+1}$ converges in $L^2$ to $2^{h+1}/(h+1)!$, and $(Q_n(h)-\E Q_n(h))/(\log n)^{h+1/2}$ has a centered normal limit with variance $D_{h,h}$, where
\begin{equation}\label{eq:diagonal}
 D_{h,h}=C_{h+1,h+1}
 =\frac{2^{2h+1}}{\bigl(h!\bigr)^2\,\bigl(2h+1\bigr)}.
\end{equation}
\end{remark}

\begin{table}[!htbp]
\centering\small
\caption{Comparison at fixed depths $r\ge1$ and $h\ge0$. The centerings use exact means.}
\label{tab:comparison}
\renewcommand{\arraystretch}{1.45}
\begin{tabular}{@{}lccc@{}}
\toprule
Quantity & Uniform $U$ & Newer $N$ & Older $O$\\
\midrule
$\E L_n(r)$ & $\sim (\log n)^r/r!$ & $\to2^r$ & $\sim2^r(\log n)^r/r!$\\
$\Var L_n(r)$ & $\displaystyle\sim\frac{(\log n)^{2r-1}}{((r-1)!)^2(2r-1)}$ & $\to K(r,r)$ & $\displaystyle\sim\frac{2^{2r-1}(\log n)^{2r-1}}{((r-1)!)^2(2r-1)}$\\
Scale for $S_n(h)$ & $(\log n)^{h+1/2}$ & $1$ & $(\log n)^{h+1/2}$\\
Fixed-level limit & Gaussian & Finite random profile & Gaussian\\
\bottomrule
\end{tabular}
\end{table}

\FloatBarrier
\section{Proofs for newer-vertex preference}\label{sec:newer}
The model superscript $N$ is suppressed in this section.

\subsection{First moments and almost-sure stabilization}\label{sec:Nfirst}
The insertion-depth formula of Lyon and Mahmoud~\cite[Theorem 6.1]{Lyon2020} is
\begin{equation}\label{eq:stirling}
 \ell_n(h)=\frac{2^h}{n!}\begin{bmatrix}n-1\\h\end{bmatrix},
 \qquad 1\le h<n,
\end{equation}
where the bracket denotes an unsigned Stirling number of the first kind. Summing over insertion times yields the expected level sizes. The two insertion-depth laws thus give the same expected level sizes in model $N$ and in the corresponding binary search tree.

For later estimates we give the weighted recurrence. The update of $W_n(h)$ and~\eqref{eq:piN} give
\begin{equation}\label{eq:FrecN}
 F_n(h)-F_{n-1}(h)=\frac{2}{n-1}F_{n-1}(h-1).
\end{equation}
Induction from $F_n(0)=1$ yields, for each fixed $h\ge0$,
\begin{equation}\label{eq:FasymN}
 F_n(h)\sim\frac{2^h}{h!}(\log n)^h.
\end{equation}
Indeed, if the result holds at level $h-1$, summing~\eqref{eq:FrecN} and using
$\sum_{j=2}^n (\log j)^{h-1}/j\sim (\log n)^h/h$ proves the next case. Consequently, for $h\ge1$,
\begin{equation}\label{eq:ellasymN}
 \ell_n(h)\sim\frac{2^h}{(h-1)!}\frac{(\log n)^{h-1}}{n^2}.
\end{equation}

The limiting means also have a short direct derivation. Put $a_h=\sum_{i\ge1}\ell_i(h)$, initially allowing the value infinity. Since $a_0=1$, Tonelli's theorem and~\eqref{eq:ellrec} give
\begin{align*}
 a_{h+1}
 &=\sum_{i\ge1}\ell_i(h)\sum_{n>i}\frac{2i}{n(n-1)}\\
 &=2a_h.
\end{align*}
Thus $a_h=2^h$. The integer-valued sequence $L_n(h)$ is nondecreasing, and monotone convergence gives $\E L_\infty(h)=2^h$. Its limit is finite almost surely. Hence the sequence is eventually constant and converges in $L^1$. Identity~\eqref{eq:identity} proves the same facts for the degree sums. The $L^2$ statement requires the next argument.

\subsection{Conditional prediction and weighted second moments}
For $h\ge1$, define
\begin{align}
 A_n(h)&=\sum_{r=0}^{h-1}2^{h-1-r}W_n(r),\label{eq:Anew}\\
 M_n(h)&=L_n(h)+\frac{2A_n(h)}n.\label{eq:Mnew}
\end{align}
Let $\mu_n(r)=\E W_n(r)$. The quantity $V_n(r,s)$ retains the definition~\eqref{eq:Vdef}.

\begin{lemma}\label{lem:weighted}
For every fixed $h\ge1$, $L_\infty(h)$ belongs to $L^2$, and
\begin{equation}\label{eq:prediction}
 M_n(h)=\E[L_\infty(h)\mid\F_n].
\end{equation}
The first and second weighted moments satisfy
\begin{align}
 \mu_n(r)&=\mu_{n-1}(r)+\frac{2}{n-1}\mu_{n-1}(r-1),\label{eq:murec}\\
 V_n(r,s)&=V_{n-1}(r,s)
 +\frac{2}{n-1}\bigl(V_{n-1}(r-1,s)+V_{n-1}(r,s-1)\bigr)\notag\\
 &\qquad+\1_{\{r=s\}}\frac{2n}{n-1}\mu_{n-1}(r-1),\label{eq:Vrec}
\end{align}
with $\mu_1(r)=\1_{\{r=0\}}$ and $V_1(r,s)=\1_{\{r=s=0\}}$. For fixed $r,s\ge1$,
\begin{equation}\label{eq:Vbound}
 V_n(r,s)=O\bigl(n\tau_n^{\min(r,s)-1}\bigr).
\end{equation}
\end{lemma}
\begin{proof}
The update $W_n(r)=W_{n-1}(r)+nI_n(r)$ proves~\eqref{eq:murec}. On multiplying two such updates, the product of their new indicators is zero unless $r=s$, in which case it is $I_n(r)$. Conditional expectation using~\eqref{eq:piN} gives~\eqref{eq:Vrec}. The first recurrence also gives
\begin{equation}\label{eq:mugen}
 \sum_{r\ge0}\mu_n(r)z^r=\prod_{j=1}^{n-1}\left(1+\frac{2z}{j}\right).
\end{equation}
It follows that $\mu_n(r)=O(\tau_n^r)$.

To prove~\eqref{eq:Vbound}, use induction on $r+s$. The boundary term is $V_n(0,s)=\mu_n(s)$. When an index in either preceding $V$ term of~\eqref{eq:Vrec} is zero, its accumulated contribution is bounded by a fixed power of $\tau_n$, which is $O(n)$. When both indices are positive, the induction hypothesis bounds each increment by $O(\tau_n^{\min(r,s)-1})$. The last term of~\eqref{eq:Vrec}, present only when $r=s$, has that same bound by~\eqref{eq:mugen}. Summing the increments proves~\eqref{eq:Vbound}. The exponent comes from the diagonal source at level $d$, where $\mu_{n-1}(d-1)=O(\tau_n^{d-1})$. A source reaching $(r,s)$ has $d\le\min(r,s)$. Propagation of a term already of order $n\tau_n^a$ divides by $n-1$ before summation and adds no logarithmic power.

For the martingale calculation, use distinct notation for the weighted insertion:
\begin{align}
 T_n(h)&=\sum_{d=1}^h2^{h-d}I_n(d),\label{eq:Tnew}\\
 \eta_n(h)&=\E[T_n(h)\mid\F_{n-1}]
            =\frac{2A_{n-1}(h)}{n(n-1)}.\label{eq:etanew}
\end{align}
The increment of $A_n(h)$ is
\[
 A_n(h)-A_{n-1}(h)
 =n\sum_{d=1}^{h-1}2^{h-1-d}I_n(d).
\]
Substitution in~\eqref{eq:Mnew} therefore gives
\begin{equation}\label{eq:martdiff}
 M_n(h)-M_{n-1}(h)=T_n(h)-\eta_n(h),
 \qquad M_1(h)=2^h.
\end{equation}
Thus $M_n(h)$ is a nonnegative martingale. Since $\E A_n(h)=O(\tau_n^{h-1})$, its difference from $L_n(h)$ tends to zero in $L^1$. Its $L^1$ limit is $L_\infty(h)$.

To improve this to $L^2$, orthogonality of martingale increments and the fact that only one depth can occur at time $n$ give
\begin{align*}
 \sum_{n=2}^\infty\E[(T_n(h)-\eta_n(h))^2]
 &\le\sum_{n=2}^\infty\E T_n(h)^2\\
 &=\sum_{d=1}^h4^{h-d}\E L_\infty(d)<\infty.
\end{align*}
The martingale consequently converges in $L^2$, and its limit agrees with the $L^1$ limit. A uniformly integrable martingale equals the conditional expectation of its terminal value, proving~\eqref{eq:prediction}.
Finally,~\eqref{eq:Vbound} gives $\E A_n(h)^2=O(n\tau_n^{h-2})$ for $h\ge2$, while $A_n(1)=1$. In both cases $2A_n(h)/n\to0$ in $L^2$. Hence $L_n(h)\to L_\infty(h)$ in $L^2$ as well.
\end{proof}

\subsection{Limiting covariance series and explicit constants}
\begin{proof}[of Theorem~\ref{thm:newer}]
The convergence and means have been proved above. For the covariance, increments of the two martingales at distinct times have zero covariance. At a common time, the conditional centering in~\eqref{eq:martdiff} gives
\begin{align}
 &\E[(T_n(h)-\eta_n(h))(T_n(k)-\eta_n(k))]\notag\\
 &\hspace{12mm}=\sum_{d=1}^{\min(h,k)}2^{h+k-2d}\ell_n(d)
                  -\E[\eta_n(h)\eta_n(k)].\label{eq:incrementcov}
\end{align}
Summing the first term over $n$ uses $\sum_n\ell_n(d)=2^d$ and gives
\[
 \sum_{d=1}^{\min(h,k)}2^{h+k-d}
 =2^{h+k}\bigl(1-2^{-\min(h,k)}\bigr).
\]
For the second term, expand both copies of $A_{n-1}$ in~\eqref{eq:etanew} and put $m=n-1$. This gives exactly the series in~\eqref{eq:Kseries}.

For positive $r,s$, estimate~\eqref{eq:Vbound} bounds its summand by a constant times $\tau_m^{\min(r,s)-1}/m^3$. Boundary terms have size $O(\tau_m^{r+s}/m^4)$. Summation of these bounds proves absolute convergence and~\eqref{eq:Ktail}. The degree-sum formula~\eqref{eq:Kdegree} follows by expanding the covariance in~\eqref{eq:identity}; deterministic terms have zero covariance.

We give the evaluation of the first constants in detail. Since the root-attachment indicators are independent,
\begin{align}
 \E W_m(1)&=2H_{m-1},\label{eq:W1mean}\\
 \Var W_m(1)&=2(m-1)+2H_{m-1}-4H_{m-1}^{(2)}.\label{eq:W1var}
\end{align}
The variance identity is obtained by summing
$j^2p^N_{j1}(1-p^N_{j1})=2j/(j-1)-4/(j-1)^2$ for $2\le j\le m$.
Therefore
\begin{align}
 A_m(1)&=1,\label{eq:A1}\\
 \E A_m(2)&=2+2H_{m-1},\label{eq:A2mean}\\
 \E A_m(2)^2&=2m+2+10H_{m-1}
               +4\bigl(H_{m-1}^2-H_{m-1}^{(2)}\bigr).\label{eq:A2second}
\end{align}
Put $d_m=m^2(m+1)^2$. The sums needed in~\eqref{eq:Kseries} are
\begin{align}
 \sum_{m\ge1}\frac1{d_m}&=2\zeta(2)-3,\label{eq:sum0}\\
 \sum_{m\ge1}\frac{H_{m-1}}{d_m}&=\zeta(2)+2\zeta(3)-4,\label{eq:sum1}\\
 \sum_{m\ge1}\frac{H_{m-1}^2-H_{m-1}^{(2)}}{d_m}
 &=2\zeta(2)+2\zeta(3)+4\zeta(4)-10,\label{eq:sum2}\\
 \sum_{m\ge1}\frac m{d_m}&=2-\zeta(2).\label{eq:summ}
\end{align}
Here are explicit reductions for these evaluations. The partial-fraction identity is
\begin{equation}\label{eq:partial}
 \frac1{d_m}=\frac1{m^2}+\frac1{(m+1)^2}-\frac2{m(m+1)}.
\end{equation}
For the harmonic numerator, shifting the second square-denominator sum and telescoping the third gives
\begin{align*}
 \sum_{m\ge1}\frac{H_{m-1}}{(m+1)^2}&=\zeta(3)+\zeta(2)-2,\\
 \sum_{m\ge1}\frac{H_{m-1}}{m(m+1)}&=1.
\end{align*}
Together with $\sum H_{m-1}/m^2=\zeta(3)$, these give~\eqref{eq:sum1}. For the quadratic numerator let $J_m=H_{m-1}^2-H_{m-1}^{(2)}$. Since $J_{m+1}-J_m=2H_{m-1}/m$, the corresponding reductions are
\begin{align*}
 \sum_{m\ge1}\frac{J_m}{m^2}&=2\zeta(4),\\
 \sum_{m\ge1}\frac{J_m}{(m+1)^2}&=2\zeta(4)+2\zeta(2)+2\zeta(3)-6,\\
 \sum_{m\ge1}\frac{J_m}{m(m+1)}&=2.
\end{align*}
The middle identity follows by shifting the index and using
\[
 \sum_{k\ge1}\frac{H_{k-1}}{k(k+1)^2}
 =3-\zeta(2)-\zeta(3).
\]
Equations~\eqref{eq:partial} and these three evaluations give~\eqref{eq:sum2}. Equations~\eqref{eq:sum0} and~\eqref{eq:summ} follow from~\eqref{eq:partial} and $1/[m(m+1)^2]=1/m-1/(m+1)-1/(m+1)^2$.

For completeness, both Euler sums used above follow from one generating function. If $e_a(m-1)$ is the elementary symmetric polynomial of degree $a$ in $1,1/2,\ldots,1/(m-1)$, then
\[
 \sum_{m\ge1}e_a(m-1)x^{m-1}
 =\frac{[-\log(1-x)]^a}{a!(1-x)}.
\]
Multiplying by $-\log x$, integrating on $(0,1)$, and substituting $y=1-x$ gives
\begin{align*}
 \sum_{m\ge1}\frac{e_a(m-1)}{m^2}
 &=\frac1{a!}\int_0^1\frac{-\log(1-y)}y(-\log y)^a\,dy\\
 &=\sum_{k\ge1}\frac1{k^{a+2}}=\zeta(a+2).
\end{align*}
All interchanges involve nonnegative terms. The cases $a=1,2$ give the required sums because $2e_2(m-1)=J_m$.

To finish the substitution, abbreviate the four sums in~\eqref{eq:sum0}--\eqref{eq:summ} as $B_0,B_1,B_2,B_m$, in their displayed order. Then
\begin{align*}
 K(1,1)&=2-4B_0,\\
 K(1,2)&=4-8B_0-8B_1,\\
 K(2,2)&=12-8B_m-8B_0-40B_1-16B_2.
\end{align*}
These expressions give~\eqref{eq:K11}--\eqref{eq:K22}. Finally,
$\Var S_\infty(1)=K(1,1)+2K(1,2)+K(2,2)$, which gives~\eqref{eq:Nsvar}.
\end{proof}

\subsection{A common limit across large depths}\label{sec:deepproof}
\begin{proof}[of Proposition~\ref{prop:large}]
For $m\ge2$ put $U_m=2^{-D(m)}$ and set $\delta_m X=X-\E[X\mid\F_{m-1}]$ for an $\F_m$-measurable square-integrable variable $X$. Divide the martingale increment~\eqref{eq:martdiff} by $2^h$. Since $T_m(h)/2^h=U_m\1_{\{D(m)\le h\}}$, its terminal sum in $L^2$ gives
\begin{equation}\label{eq:deepmart}
 2^{-h}L_\infty(h)=1+\sum_{m\ge2}\delta_m\!\left(U_m\1_{\{D(m)\le h\}}\right).
\end{equation}
By Tonelli's theorem and $\E L_\infty(d)=2^d$,
\[
 \sum_{m\ge2}\E U_m^2
 =\sum_{d\ge1}4^{-d}\E L_\infty(d)=1.
\]
The centered increments are orthogonal at distinct times and conditional centering contracts the $L^2$ norm. Thus $Y=1+\sum_{m\ge2}\delta_m U_m$ exists in $L^2$, has mean one and variance $\kappa=\sum_{m\ge2}\E\Var(U_m\mid\F_{m-1})$. Subtracting~\eqref{eq:deepmart} from this series gives
\[
 \E\big|Y-2^{-h}L_\infty(h)\big|^2
 \le\sum_{m\ge2}\E[U_m^2\1_{\{D(m)>h\}}]
 =\sum_{d>h}2^{-d}=2^{-h}.
\]
To evaluate $\kappa$, put $P_m(z)=\sum_{d\ge0}W_m(d)z^d$,
$a_m=\E[P_m(1/2)^2]$, $b_m=a_m/[m(m+1)]$ and $g_m=\E P_m(1/4)$.
The parent rule gives $\E[U_{m+1}\mid\F_m]=P_m(1/2)/[m(m+1)]$.
Subtracting the squared conditional means from $\sum_{m\ge2}\E U_m^2=1$ gives
\[
 \kappa=1-\sum_{m\ge1}\frac{b_m}{m(m+1)}.
\]
Equations~\eqref{eq:Vrec} and~\eqref{eq:mugen} give $b_1=1/2$, $g_1=1$ and
\[
 b_{m+1}-b_m=\dfrac{g_m}{2m(m+2)},\qquad
 g_{m+1}=\left(1+\dfrac{1}{2m}\right)g_m .
\]
The binomial series gives $G(t)=\sum_{m\ge1}g_m t^{m-1}=(1-t)^{-3/2}$ for $0\le t<1$.
Because the terms are nonnegative, summing the increments of $b_m$ and integrating term by term yields
\[
 b_\infty=\frac12+\frac14\int_0^1(1-t^2)G(t)\,dt
         =\frac12+\frac14\int_0^1\frac{1+t}{\sqrt{1-t}}\,dt=\frac43 .
\]
In particular $b_m$ is bounded. Summation by parts, or Tonelli applied to the increments of $b_m$, now gives
\begin{align*}
 \sum_{m\ge1}\frac{b_m}{m(m+1)}
 &= b_1+\sum_{m\ge1}\frac{b_{m+1}-b_m}{m+1}
  = \frac12+\frac12\sum_{m\ge1}\frac{g_m}{m(m+1)(m+2)}\\
 &=\frac12+\frac14\int_0^1(1-t)^2G(t)\,dt=\frac23 .
\end{align*}
The integral uses $\int_0^1t^{m-1}(1-t)^2\,dt=2/[m(m+1)(m+2)]$.
Hence $\kappa=1/3$ and $\E Y^2=4/3$. Covariance convergence in~\eqref{eq:deepK} follows from~\eqref{eq:deepL} by Cauchy--Schwarz.
\end{proof}

\paragraph{A one-sided error bound.}
Write $q_m=\E[U_m\mid\F_{m-1}]$ and $q_{m,h}=\E[U_m\1_{\{D(m)\le h\}}\mid\F_{m-1}]$.
Using~\eqref{eq:deepmart} for levels $h$ and $k$ and subtracting their covariance from $\kappa$ gives, for $d=\min(h,k)$,
\[
 \frac13-\frac{K(h,k)}{2^{h+k}}
 =2^{-d}-\sum_{m\ge2}\E\left[q_m^2-q_{m,h}q_{m,k}\right]\le 2^{-d}.
\]
Indeed the uncentered product omits precisely the depths above $d$, and $0\le q_{m,h},q_{m,k}\le q_m$. This bounds the shortfall below $1/3$. The sharp two-sided rate remains open here.

\paragraph{Support of the limit.}
The martingale $P_m(1/2)/m$ starts at one. At step $m+1$ its increment is $U_{m+1}-P_m(1/2)/[m(m+1)]$, so $P_m(1/2)/m=\E[Y\mid\F_m]$.
On the positive-probability event that vertices $2,\ldots,m$ all attach to the root, this conditional mean is $(m+1)/4+1/(2m)$. It can exceed any fixed threshold. Consequently $Y$ has unbounded support, which rules out the uniform law on $[0,2]$ despite the agreement of its first two moments.

\paragraph{Three further covariance entries.}
At cutoff $M=2^{20}$, the interval computation described below gives
\begin{equation}\label{eq:enclosures}
\begin{aligned}
 2.39112 &\le K(1,3)\le 2.39115,\\
 7.89110 &\le K(2,3)\le 7.89115,\\
 19.88516&\le K(3,3)\le19.88525.
\end{aligned}
\end{equation}
To certify the omitted positive correction in~\eqref{eq:Kseries}, note that $2^{-h}A_m(h)\le P_m(1/2)/2$ and $b_m\le4/3$. Thus its tail after $M$ is at most
\[
 2^{h+k}\sum_{m>M}\frac{b_m}{m(m+1)}
 \le \frac{4\,2^{h+k}}{3(M+1)} .
\]
For the first $M$ terms the companion script propagates lower and upper integer bounds for the nonnegative recurrences~\eqref{eq:murec}--\eqref{eq:Vrec}, rounding outwards on a grid of spacing $10^{-30}$. It subtracts the accumulated correction interval and the one-sided tail bound before rounding the displayed endpoints outwards. The proof of the enclosures uses only the finite recurrences and the tail estimate. The six values in~\eqref{eq:K11}--\eqref{eq:K22} and~\eqref{eq:enclosures} all satisfy $2^{-h-k}K(h,k)<1/3$. Does $K(h,k)\le 2^{h+k}/3$ hold for every $h,k\ge1$?

\subsection{Stabilization rates}
\begin{proof}[of Theorem~\ref{thm:tails}]
Let
\begin{equation}\label{eq:tailmean}
 b_n(d)=\E[L_\infty(d)-L_n(d)].
\end{equation}
By~\eqref{eq:ellasymN} and comparison of the summable tail with an integral,
\begin{equation}\label{eq:btail}
 b_n(d)\sim\frac{2^d}{(d-1)!}\frac{(\log n)^{d-1}}n,
 \qquad d\ge1.
\end{equation}
For $d=1$, the exact value is $b_n(1)=2/n$.

Write $\Delta_n(h)=L_\infty(h)-L_n(h)$. Conditional expectation~\eqref{eq:prediction} makes $L_\infty(h)-M_n(h)$ orthogonal in $L^2$ to the $\F_n$-measurable variable $M_n(h)-L_n(h)$. Hence
\begin{align}
 \E\Delta_n(h)^2
 &=\E[(L_\infty(h)-M_n(h))^2]+\frac4{n^2}\E A_n(h)^2\notag\\
 &=\sum_{d=1}^h4^{h-d}b_n(d)
   -\sum_{m>n}\E\eta_m(h)^2+\frac4{n^2}\E A_n(h)^2.\label{eq:tailorth}
\end{align}
The second equality follows by summing the squared martingale increments in~\eqref{eq:martdiff} after time $n$.
For $h\ge2$, Lemma~\ref{lem:weighted} gives
\begin{align}
 \E A_n(h)^2&=O(n\tau_n^{h-2}),\label{eq:Abound}\\
 \sum_{m>n}\E\eta_m(h)^2&=O(\tau_n^{h-2}/n^2).\label{eq:etabound}
\end{align}
Each term $b_n(d)$ with $d<h$ is $o(b_n(h))$. The last two terms of~\eqref{eq:tailorth} are also $o(b_n(h))$, by~\eqref{eq:btail}--\eqref{eq:etabound}. Thus $\E\Delta_n(h)^2\sim b_n(h)$. When $h=1$, $A_n(1)=1$ and the two correction terms have orders $n^{-3}$ and $n^{-2}$, so the same conclusion holds.

For a nonnegative integer-valued random variable $Z$,
\begin{equation}\label{eq:integerineq}
 0\le\E Z-\Pp(Z>0)\le\tfrac12\E[Z(Z-1)].
\end{equation}
Apply this to $\Delta_n(h)$. Since $\E\Delta_n(h)^2\sim\E\Delta_n(h)$, the right-hand side is $o(\E\Delta_n(h))$. This proves~\eqref{eq:jointtail} for the level count.

For degree sums,~\eqref{eq:identity} gives
\[
 S_\infty(h)-S_n(h)=\Delta_n(h)+\Delta_n(h+1),
\]
where $\Delta_n(0)=0$. For $h\ge1$, the ratio of the second moments of these two terms tends to zero. More explicitly, Cauchy--Schwarz gives
\[
 \frac{\E[\Delta_n(h)\Delta_n(h+1)]}{\E\Delta_n(h+1)^2}
 \le\left(\frac{\E\Delta_n(h)^2}{\E\Delta_n(h+1)^2}\right)^{1/2}
 =O((\log n)^{-1/2}).
\]
The square of the sum and its mean are therefore both asymptotic to $b_n(h+1)$. Applying~\eqref{eq:integerineq} to that sum proves~\eqref{eq:jointtail} for the degree sum. At $h=0$ it equals $\Delta_n(1)$. For either increment, $2\Pp(\Delta\ge2)\le\E[\Delta(\Delta-1)]=o(b_n(d))$ and $\Pp(\Delta>0)\sim b_n(d)$. Their ratio gives the conditional one-increment assertion.
\end{proof}

\section{Proofs for older-vertex preference}\label{sec:older}
The superscript $O$ is suppressed inside the proofs. We use $X_n$ for centered level counts, $Z_n$ for centered weighted counts and $\alpha_n$ for iterated root noise; these are distinct from the prediction variables in Section~\ref{sec:newer}.

\subsection{First moments}
\begin{lemma}\label{lem:Ofirst}
For fixed $r\ge0$, $\E L_n(r)=O(\tau_n^r)$. For fixed $r\ge1$,
\begin{align}
 F_n(r)&=O(n\tau_n^{r-1}),\label{eq:OFbound}\\
 \ell_n(r)&\sim\frac{2^r}{(r-1)!}\frac{(\log n)^{r-1}}n.\label{eq:Oellasym}
\end{align}
The mean assertions~\eqref{eq:Omean} and~\eqref{eq:Osmean} follow.
\end{lemma}
\begin{proof}
At the first level,
\begin{align*}
 \ell_n(1)&=2/n,\\
 \E L_n(1)&=2(H_n-1),\\
 F_n(1)&=2(n-1).
\end{align*}
Taking expectations in~\eqref{eq:piO} gives
\begin{equation}\label{eq:Oellrec}
 \ell_n(r)=\frac{2\E L_{n-1}(r-1)}{n-1}
            -\frac{2F_{n-1}(r-1)}{n(n-1)}.
\end{equation}
Assume the level-count bound at $r-1$. Nonnegativity of the second term implies $\ell_n(r)=O(\tau_n^{r-1}/n)$. Summing this bound over $n$ proves the next level-count bound. Summing it after multiplication by $n$ proves~\eqref{eq:OFbound}.

We now prove the leading constants by induction. At $r=1$ they are explicit above. If~\eqref{eq:Omean} holds at $r-1\ge1$, the first term of~\eqref{eq:Oellrec} is asymptotic to $2^r(\log n)^{r-1}/((r-1)!n)$. Estimate~\eqref{eq:OFbound} makes the second term $O(\tau_n^{r-2}/n)$, a smaller order. This proves~\eqref{eq:Oellasym}; summation proves~\eqref{eq:Omean} at level $r$. In~\eqref{eq:identity}, the mean of level $h+1$ exceeds that of level $h$ by a logarithmic order. This yields~\eqref{eq:Osmean}, including the root case.
\end{proof}

\subsection{Centered recurrences and moment bounds}
Define
\begin{align}
 X_n(r)&=L_n(r)-\E L_n(r),\label{eq:X}\\
 Z_n(r)&=W_n(r)-\E W_n(r).\label{eq:Z}
\end{align}
Both quantities are zero at level zero. With the conditional probability $\pi_n^O(r)$ from~\eqref{eq:piO}, put
\begin{equation}\label{eq:eps}
 \varepsilon_n(r)=I_n(r)-\pi_n^O(r).
\end{equation}
For each fixed $r$, these are martingale differences. Hence increments at different times are orthogonal in $L^2$, and
\begin{equation}\label{eq:epsbound}
 \E\varepsilon_n(r)^2
 =\E[\pi_n^O(r)(1-\pi_n^O(r))]\le\ell_n(r).
\end{equation}
Subtracting the expectations from $L_n(r)=L_{n-1}(r)+I_n(r)$ and $W_n(r)=W_{n-1}(r)+nI_n(r)$ gives
\begin{align}
 X_n(r)-X_{n-1}(r)
 &=\frac2{n-1}X_{n-1}(r-1)
   -\frac2{n(n-1)}Z_{n-1}(r-1)+\varepsilon_n(r),\label{eq:Xrec}\\
 Z_n(r)-Z_{n-1}(r)
 &=\frac{2n}{n-1}X_{n-1}(r-1)
   -\frac2{n-1}Z_{n-1}(r-1)+n\varepsilon_n(r).\label{eq:Zrec}
\end{align}

\begin{lemma}\label{lem:moments}
For each fixed $r\ge1$,
\begin{equation}\label{eq:Xbound}
 \norm{X_n(r)}=O(\tau_n^{r-1/2}).
\end{equation}
The weighted counts satisfy
\begin{align}
 \norm{Z_n(1)}&=O(n),\label{eq:Zrootbound}\\
 \norm{Z_n(r)}&=O(n\tau_n^{r-3/2}),\qquad r\ge2.\label{eq:Zbound}
\end{align}
\end{lemma}
\begin{proof}
Orthogonality,~\eqref{eq:epsbound} and Lemma~\ref{lem:Ofirst} give two estimates for the sums of martingale differences:
\begin{align}
 \left\lVert\sum_{j=2}^n\varepsilon_j(r)\right\rVert_2^2
 &\le\sum_{j=2}^n\ell_j(r)=O(\tau_n^r),\label{eq:unweightednoise}\\
 \left\lVert\sum_{j=2}^nj\varepsilon_j(r)\right\rVert_2^2
 &\le\sum_{j=2}^nj^2\ell_j(r)=O(n^2\tau_n^{r-1}).\label{eq:weightednoise}
\end{align}
When $r=1$, both drift terms in~\eqref{eq:Xrec} and~\eqref{eq:Zrec} vanish, so these estimates prove the initial cases.

Suppose the assertions hold at level $r-1$, where $r\ge2$. Sum~\eqref{eq:Xrec} and display its three terms:
\begin{equation}\label{eq:Xthree}
 X_n(r)=2\sum_{j=2}^n\frac{X_{j-1}(r-1)}{j-1}
       -2\sum_{j=2}^n\frac{Z_{j-1}(r-1)}{j(j-1)}
       +\sum_{j=2}^n\varepsilon_j(r).
\end{equation}
By the induction hypothesis and Minkowski's inequality, the norm of the first sum, including its factor $2$, is bounded by
\begin{equation}\label{eq:Xdrift1}
 O\left(\sum_{j=2}^n\frac{\tau_j^{r-3/2}}j\right)
 =O(\tau_n^{r-1/2}).
\end{equation}
The norm of the second term is
\begin{equation}\label{eq:Xdrift2}
 \begin{cases}
 O\bigl(\sum_{j=2}^n j^{-1}\bigr)=O(\tau_n),&r=2,\\
 O\bigl(\sum_{j=2}^n\tau_j^{r-5/2}/j\bigr)
       =O(\tau_n^{r-3/2}),&r\ge3.
 \end{cases}
\end{equation}
The norm of the third term is $O(\tau_n^{r/2})$ by~\eqref{eq:unweightednoise}. Every bound is $O(\tau_n^{r-1/2})$, proving~\eqref{eq:Xbound} at level $r$.

For the weighted count, the corresponding decomposition is
\begin{equation}\label{eq:Zthree}
 Z_n(r)=2\sum_{j=2}^n\frac{jX_{j-1}(r-1)}{j-1}
       -2\sum_{j=2}^n\frac{Z_{j-1}(r-1)}{j-1}
       +\sum_{j=2}^nj\varepsilon_j(r).
\end{equation}
The first term has norm
\begin{equation}\label{eq:Zdrift1}
 O\left(\sum_{j=2}^n\tau_j^{r-3/2}\right)
 =O(n\tau_n^{r-3/2}).
\end{equation}
The second term has norm $O(n)$ for $r=2$ and $O(n\tau_n^{r-5/2})$ for $r\ge3$. By~\eqref{eq:weightednoise}, the third term has norm $O(n\tau_n^{(r-1)/2})$. For $r\ge2$, $(r-1)/2\le r-3/2$. All three terms in~\eqref{eq:Zthree} are therefore $O(n\tau_n^{r-3/2})$. This closes the induction.
\end{proof}

\subsection{Approximation by root-attachment noise}
The root-attachment probability is $2/j$. Define
\begin{align}
 \xi_j&=\1_{\{\Gamma_j=1\}}-\frac2j,\label{eq:xi}\\
 v_j&=\E\xi_j^2=\frac2j-\frac4{j^2}.\label{eq:vj}
\end{align}
The $\xi_j$ are independent, centered, and satisfy $|\xi_j|\le1$. In particular,
$X_n(1)=\sum_{j=2}^n\xi_j$. Let
\begin{align}
 \alpha_n(1)&=\sum_{j=2}^n\xi_j,\label{eq:alpha1}\\
 \alpha_n(r)&=2\sum_{j=2}^n\frac{\alpha_{j-1}(r-1)}{j-1},\qquad r\ge2.\label{eq:alpharec}
\end{align}
Empty sums are zero.

\begin{lemma}\label{lem:rootapprox}
For every fixed $r\ge2$,
\begin{equation}\label{eq:alphaerror}
 \norm{X_n(r)-\alpha_n(r)}=O(\tau_n^{r-1}).
\end{equation}
\end{lemma}
The fluctuation scale at level $r$ is $(\log n)^{r-1/2}$. Dividing the error bound by this scale gives $O(1/\sqrt{\log n})$. This slow decay is consistent with the finite-size variance discrepancies in Section~\ref{sec:simulation} (Figure~\ref{fig:older} and Table~\ref{tab:finiteO}). The estimate controls the mean-square remainder in the representation.

\begin{proof}
Write $\beta_n(r)=X_n(r)-\alpha_n(r)$, so $\beta_n(1)=0$. Subtracting~\eqref{eq:alpharec} from~\eqref{eq:Xthree} gives the exact decomposition
\begin{equation}\label{eq:betathree}
 \beta_n(r)=2\sum_{j=2}^n\frac{\beta_{j-1}(r-1)}{j-1}
       -2\sum_{j=2}^n\frac{Z_{j-1}(r-1)}{j(j-1)}
       +\sum_{j=2}^n\varepsilon_j(r).
\end{equation}
For $r=2$, the first term is zero, the second has norm $O(\tau_n)$ by~\eqref{eq:Zrootbound}, and the third has norm $O(\tau_n)$ by~\eqref{eq:unweightednoise}. This proves the initial case.

For $r\ge3$, assume~\eqref{eq:alphaerror} at $r-1$. The norm of the first term of~\eqref{eq:betathree} is bounded by
\[
 O\left(\sum_{j=2}^n\frac{\tau_j^{r-2}}j\right)=O(\tau_n^{r-1}).
\]
The second term has norm $O(\tau_n^{r-3/2})$ by~\eqref{eq:Zbound}; the third has norm $O(\tau_n^{r/2})$. Since $r/2\le r-1$, their sum is $O(\tau_n^{r-1})$. Minkowski's inequality now proves the induction step.
\end{proof}

Iterating~\eqref{eq:alpharec} gives
\begin{equation}\label{eq:alphakernel}
 \alpha_n(r)=2^{r-1}\sum_{j=2}^n\xi_j
 e_{r-1}\left(\frac1j,\frac1{j+1},\ldots,\frac1{n-1}\right).
\end{equation}
Here $e_d$ is the elementary symmetric polynomial of degree $d$, $e_0=1$, and $e_d=0$ when fewer than $d$ arguments are present. To check~\eqref{eq:alphakernel} inductively, substitute its expression at level $r-1$ into~\eqref{eq:alpharec}. Each contribution from $\xi_j$ is multiplied by the reciprocal of a new, larger index. Thus every increasing $(r-1)$-tuple from $\{j,\ldots,n-1\}$ occurs once.

\begin{lemma}\label{lem:kernel}
For fixed $d\ge1$, uniformly over $2\le j\le n$,
\begin{equation}\label{eq:eapprox}
 e_d\left(\frac1j,\ldots,\frac1{n-1}\right)
 =\frac{\log^d(n/j)}{d!}
  +O\left(\frac{(1+\log(n/j))^{d-1}}j\right).
\end{equation}
\end{lemma}
\begin{proof}
Put $u=\sum_{k=j}^{n-1}k^{-1}$. Integral comparison gives
$u=\log(n/j)+O(1/j)$, uniformly in $j$. In the expansion of $u^d$, all ordered tuples of distinct indices contribute exactly $d!e_d$. If $d\ge2$, every other tuple contains at least one equal pair of indices. Choosing that pair and summing the remaining indices bounds its total contribution by
\[
 0\le u^d-d!e_d
 \le \binom{d}{2}u^{d-2}\sum_{k=j}^{n-1}\frac1{k^2}.
\]
The square-reciprocal sum is $O(1/j)$. Also,
$|u^d-\log^d(n/j)|=O((1+\log(n/j))^{d-1}/j)$ by the difference-of-powers identity. These two bounds give~\eqref{eq:eapprox}. For $d=1$, the harmonic-sum estimate alone suffices. Empty intervals satisfy the same estimate directly.
\end{proof}

By independence and~\eqref{eq:vj}, the squared $L^2$ error from replacing the elementary symmetric polynomial in~\eqref{eq:alphakernel} by its logarithmic term is, for $r\ge2$, at most
\begin{equation}\label{eq:kernelL2}
 O\left(\sum_{j=2}^n\frac{(1+\log(n/j))^{2r-4}}{j^3}\right)
 =O(\tau_n^{2r-4}).
\end{equation}
Combining this with Lemma~\ref{lem:rootapprox}, we obtain
\begin{equation}\label{eq:rootrepresentation}
 X_n(r)=P_n(r)+R_n(r),
\end{equation}
where
\begin{align}
 P_n(r)&=\frac{2^{r-1}}{(r-1)!}\sum_{j=2}^n\xi_j\log^{r-1}(n/j),\label{eq:P}\\
 \norm{R_n(r)}&=O(\tau_n^{r-1}),\qquad r\ge2.\label{eq:Rbound}
\end{align}
At $r=1$, $R_n(1)=0$, and the zeroth power in~\eqref{eq:P} is interpreted as $1$, including at $j=n$.

\subsection{Covariance constants and the joint limit}
\begin{proof}[of Theorem~\ref{thm:older}]
The mean assertions were proved in Lemma~\ref{lem:Ofirst}. For an integer $q\ge0$, we first record
\begin{equation}\label{eq:logsum}
 \sum_{j=2}^nv_j\log^q(n/j)
 =\frac2{q+1}(\log n)^{q+1}+O(\tau_n^q).
\end{equation}
For $q=0$, every zeroth power is $1$. The function $x^{-1}\log^q(n/x)$ is nonnegative and nonincreasing on $[1,n]$. Consequently its sum--integral error, including the omitted term at $j=1$, is $O(\tau_n^q)$. Its integral is
\[
 \int_1^n\frac{\log^q(n/x)}x\,dx=\frac{(\log n)^{q+1}}{q+1}.
\]
The part from $4/j^2$ in~\eqref{eq:vj} is bounded by
$4(\log n)^q\sum_{j\ge2}j^{-2}=O(\tau_n^q)$. These estimates prove~\eqref{eq:logsum}.

Independence of the $\xi_j$ and~\eqref{eq:P} give
\begin{equation}\label{eq:Pcov}
 \Cov(P_n(r),P_n(s))
 =\frac{2^{r+s-2}}{(r-1)!(s-1)!}
   \sum_{j=2}^nv_j\log^{r+s-2}(n/j).
\end{equation}
Equation~\eqref{eq:logsum} identifies its leading term as
$C_{r,s}(\log n)^{r+s-1}$. The remainders do not change that term. For example, when $s\ge2$,
\[
 |\E[P_n(r)R_n(s)]|
 \le\norm{P_n(r)}\norm{R_n(s)}
 =O(\tau_n^{r+s-3/2})=o((\log n)^{r+s-1}).
\]
The other cross term satisfies the same bound, and the product of two nonzero remainders is $O(\tau_n^{r+s-2})$. Remainders at level $1$ vanish. This proves~\eqref{eq:Ocov}.

For the joint limit, fix real coefficients $a_1,\ldots,a_m$. The linear combination of the normalized $P_n(r)$ is a triangular-array sum $\sum_{j=2}^nb_{n,j}\xi_j$ with
\begin{equation}\label{eq:bnj}
 b_{n,j}=\frac1{\sqrt{\log n}}
   \sum_{r=1}^m\frac{a_r2^{r-1}}{(r-1)!}
   \left(\frac{\log(n/j)}{\log n}\right)^{r-1}.
\end{equation}
Since the fraction inside the powers lies in $[0,1]$,
$\max_j|b_{n,j}|=O((\log n)^{-1/2})$. The bound $|\xi_j|\le1$ then makes every Lindeberg truncation term zero for all sufficiently large $n$. By~\eqref{eq:logsum}, the sum of variances converges to $\sum_{r,s}a_ra_sC_{r,s}$.

We explicitly verify that this limit is positive for a nonzero coefficient vector. Define
\[
 f_r(u)=\frac{2^{r-1}}{(r-1)!}u^{r-1}.
\]
From~\eqref{eq:C} and $\int_0^1u^{r+s-2}\,du=1/(r+s-1)$,
\begin{equation}\label{eq:gramentry}
 C_{r,s}=2\int_0^1f_r(u)f_s(u)\,du.
\end{equation}
Multiplying by $a_ra_s$ and summing yields
\begin{equation}\label{eq:gram}
 \sum_{r,s=1}^ma_ra_sC_{r,s}
 =2\int_0^1\left(\sum_{r=1}^ma_rf_r(u)\right)^2\,du.
\end{equation}
Equivalently $C=2\Lambda H\Lambda$ on any finite block, where $\Lambda_{rr}=2^{r-1}/(r-1)!$ and $H_{rs}=1/(r+s-1)$. This is the Hilbert matrix in the uniform covariance~\eqref{eq:ucov}, the Gram matrix of $1,u,\ldots,u^{m-1}$ on $[0,1]$. The polynomials $f_r$ have distinct degrees. Their linear combination is identically zero only when all $a_r$ are zero; otherwise its squared integral is positive. This establishes positive definiteness and gives the origin of the integral representation.

The Lindeberg--Feller theorem now gives the normal limit for every nonzero linear combination. The Cram\'er--Wold device gives the multivariate limit for $P_n$, and~\eqref{eq:Rbound} gives
\[
 \frac{\norm{R_n(r)}}{(\log n)^{r-1/2}}=O((\log n)^{-1/2})\longrightarrow0.
\]
Slutsky's theorem proves~\eqref{eq:Oclt}.

Finally, the centered version of~\eqref{eq:identity} is
\[
 S_n(h)-\E S_n(h)=X_n(h)+X_n(h+1).
\]
For $h\ge1$, the first term has norm $O(\tau_n^{h-1/2})$ and the second has fluctuation scale $(\log n)^{h+1/2}$. At $h=0$ the first term is zero. Thus the leading degree-sum covariance is $C_{h+1,k+1}=D_{h,k}$, and its joint Gaussian limit follows from the corresponding subvector in~\eqref{eq:Oclt}. Positive definiteness follows either from~\eqref{eq:gram} or from the same monomial argument. This proves~\eqref{eq:Oscov}--\eqref{eq:Osvar}. Dividing the variance estimates by the squares of the leading means gives a quantity of order $1/\log n$; together with the mean limits this proves~\eqref{eq:OL2}--\eqref{eq:OSL2}.
\end{proof}

\section{Numerical illustrations}\label{sec:simulation}
For each model we generated $10{,}000$ independent growth paths. On each path we recorded the first three nonzero levels at $n=10^2,10^3,10^4,10^5$. The simulations use the full parent-choice construction. For model $N$, inverse sampling uses the exact cumulative probability
\begin{equation}\label{eq:cdf}
 \Pp(\Gamma_n\le i)=\frac{i(i+1)}{n(n-1)},\qquad 1\le i<n.
\end{equation}
Reflecting this draw as $i\mapsto n-i$ gives model $O$. The implementation uses \texttt{std::mt19937\_64}, with seeds $2026111301$ and $2026111302$, respectively. Floating-point inversion is followed by a check against the adjacent cumulative-probability boundaries. Vertex depths larger than three can be capped at four because a descendant of such a vertex cannot return to one of the recorded levels.

Exact finite-$n$ means and second moments were computed separately from the recurrence in Appendix~\ref{app:momentrec}. To verify that recurrence, we compared it with exhaustive enumeration of all parent sequences for $2\le n\le7$, for $N$, $O$ and $U$. The largest absolute discrepancy in a state mean or covariance was below $6\cdot10^{-16}$ in the numerical calculation. The simulation means at all four recorded sizes were within $1.70$ theoretical Monte Carlo standard errors of their exact values. 

\begin{figure}[!htbp]
 \centering
 \includegraphics[width=\textwidth]{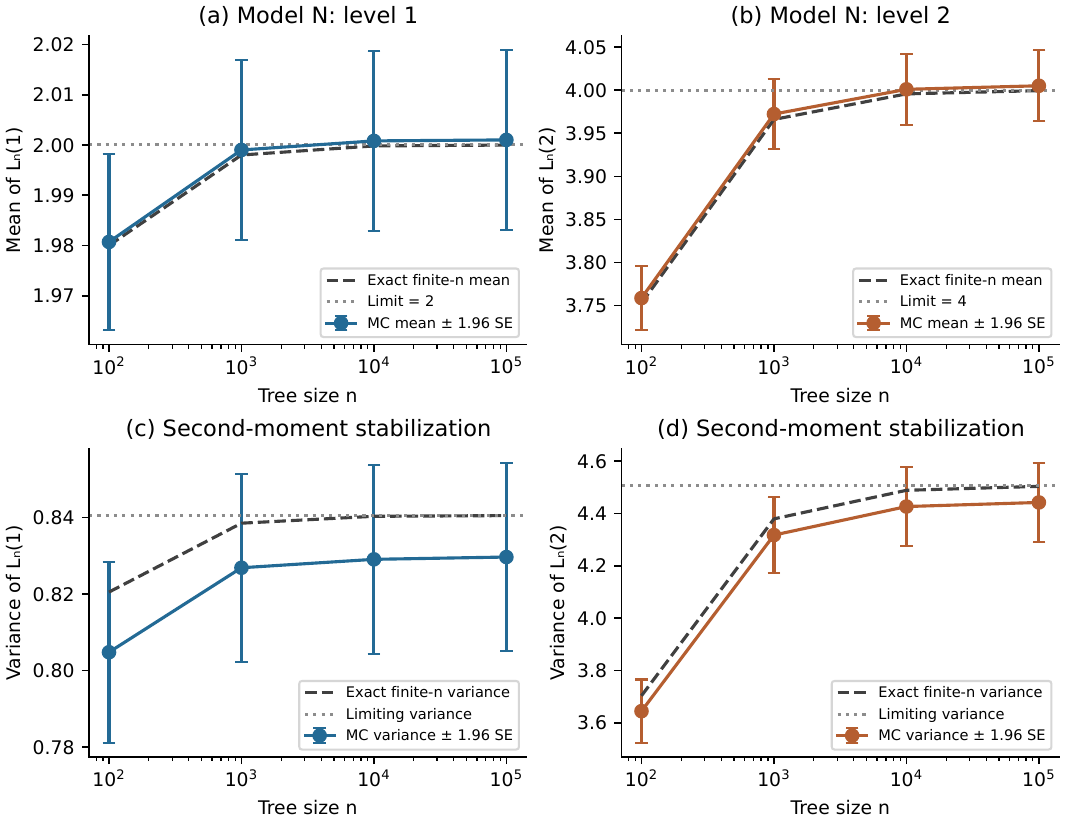}
 \caption{Model $N$: empirical means and variances of levels $1$ and $2$ from $10{,}000$ independent growth paths. Dashed curves are exact finite-$n$ moments computed by recurrence; dotted lines are the limiting constants. Error bars are $1.96$ estimated Monte Carlo standard errors, using empirical fourth moments for the variance panels. All four sizes use the same paths. In (c,d), sample variances use the unbiased estimator across paths. Estimates at successive sizes are correlated because most paths have already stabilized, so several discrepancies with the same sign can reflect the same sampling error.}\label{fig:newer}
\end{figure}

Figure~\ref{fig:newer} shows the mean and variance stabilization of $L_n^N(1)$ and $L_n^N(2)$. The exact finite-size curves make the approach to the constants visible even where Monte Carlo uncertainty is comparable to the remaining bias. Level zero is omitted because its count is identically one. The root degree is nevertheless represented, since $S_n^N(0)=L_n^N(1)$.

\begin{figure}[!htbp]
 \centering
 \includegraphics[width=\textwidth]{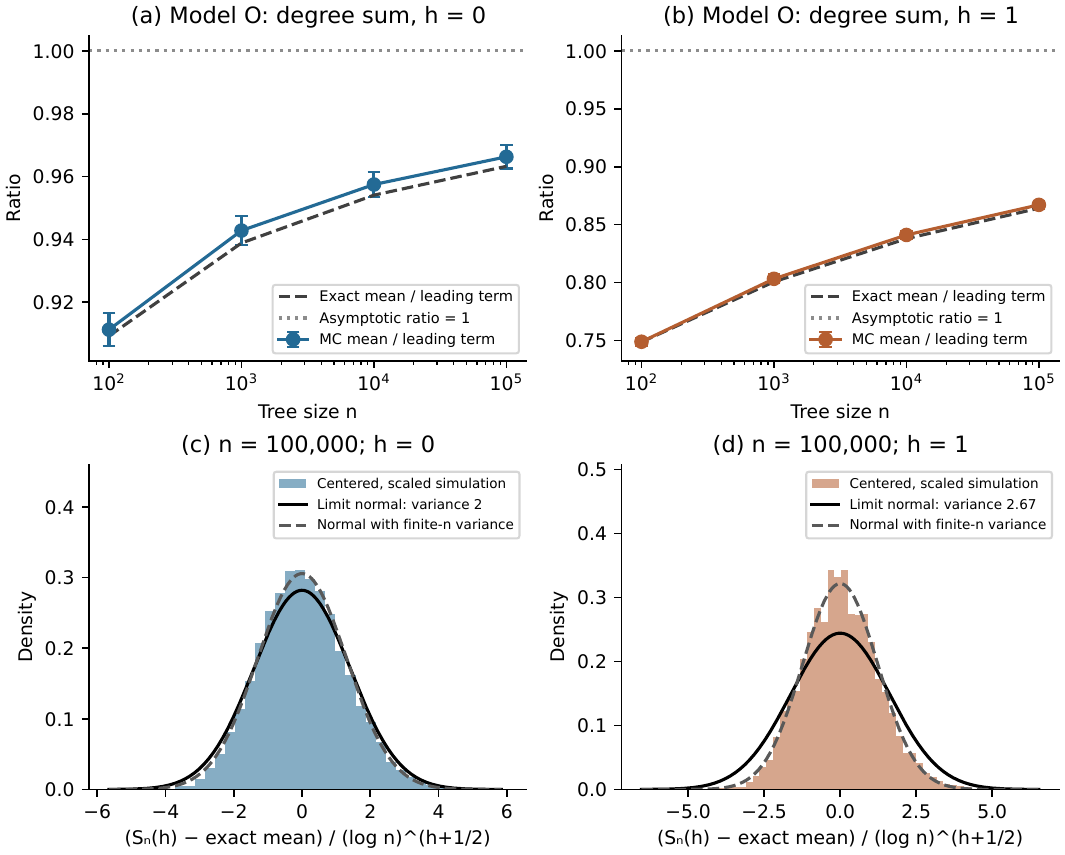}
 \caption{Model $O$: ratios of empirical and exact means to the leading terms, and centered, logarithmically normalized degree sums. Histograms use $10{,}000$ independent trees of size $10^5$. The horizontal references in (a,b) are at ratio $1$. In (c,d) the solid Gaussian curves have variances $D_{0,0}=2$ and $D_{1,1}=8/3\approx2.67$, respectively. The dashed Gaussian curves use the exact finite-size variances on the same scale.}\label{fig:older}
\end{figure}

For model $O$, the upper panels of Figure~\ref{fig:older} compare the means of $S_n^O(0)$ and $S_n^O(1)$ with their leading asymptotic terms. The lower panels use exact finite-$n$ centering and the normalization $(\log n)^{h+1/2}$ from Theorem~\ref{thm:older}. They include both the limiting Gaussian density and a Gaussian density with the exact finite-$n$ variance at that normalization.

The distinction between these two curves is substantial. At $n=10^5$, the exact variance divided by its leading asymptotic term is approximately $0.8512$ for $h=0$ and $0.5762$ for $h=1$. Thus a visibly narrower histogram at this size is compatible with the asymptotic theorem: the variance itself has not reached its leading approximation. Table~\ref{tab:finiteO} reports the underlying values. The finite-variance Gaussian curve isolates the effect of this variance correction.

\begin{table}[!htbp]
\centering\small
\caption{Degree sums in model $O$ at $n=100{,}000$. MC denotes the $10{,}000$-replicate Monte Carlo calculation; parentheses give the MC standard error of the mean. The last column is exact variance divided by its asymptotic leading term.}
\label{tab:finiteO}
\begin{tabular}{@{}cccccc@{}}
\toprule
$h$ & Exact mean & MC mean (SE) & Exact variance & MC variance & Variance ratio\\
\midrule
0 & 22.1803 & 22.2501 (0.0437) & 19.6006 & 19.0941 & 0.8512\\
1 & 229.0918 & 229.8584 (0.4853) & 2,344.5763 & 2,355.5045 & 0.5762\\
\bottomrule
\end{tabular}
\end{table}

\FloatBarrier

\section{Concluding remarks}\label{sec:concl}
The covariance kernel in model $N$ has two complementary descriptions. Its low entries distinguish the first few levels, while its large-depth asymptotic shows that normalized level counts share a single random limit. The common limit has variance $1/3$ and unbounded support. Its full distribution and the existence of zeta-value formulas for all finite entries remain open here. In model $O$ the normalized profile has a Gaussian limit driven by attachments at the root. At $10^5$ vertices the exact variance of the degree sum $S_n^O(1)$ is still only about $58\%$ of its leading approximation.

The large-depth limit in Proposition~\ref{prop:large} concerns the completed infinite tree. Passing to a depth that increases along with the size of a finite tree would require a stabilization estimate uniform over that depth range. Other natural extensions include power weights~\cite{Lyon2022} and age weights $(n-i)^\alpha$.

\appendix
\renewcommand*{\theHsection}{appendix.\Alph{section}}
\renewcommand*{\theHequation}{\theHsection.\arabic{equation}}
\section{Finite-size identities and computational checks}\label{app:checks}
\subsection{The first two levels}
The following formulas fix the centering and conventions used in Section~\ref{sec:simulation}. They follow from the marginal recurrences and, for model $N$, also from~\cite[Theorem 6.1]{Lyon2020}. For $n\ge2$,
\begin{align}
 \E L_n^N(1)&=2\left(1-\frac1n\right),\label{eq:finiteN1}\\
 \E L_n^N(2)&=4\left(1-\frac1{n-1}-\frac{H_{n-2}}n\right),\label{eq:finiteN2}\\
 \E L_n^O(1)&=2(H_n-1),\label{eq:finiteO1}\\
 \E L_n^O(2)&=2H_n^2+2H_n^{(2)}-8H_n-\frac{4H_n}n+8.\label{eq:finiteO2}
\end{align}
Identity~\eqref{eq:identity} gives $\E S_n(0)=\E L_n(1)$ and $\E S_n(1)=\E L_n(1)+\E L_n(2)$.

For the newer-vertex second level,~\eqref{eq:ellrec} gives
$\ell_n^N(2)=4H_{n-2}/[n(n-1)]$. The telescoping step in its summation is, with $m=n-2$,
\begin{align*}
 \sum_{i=1}^{m}\frac{H_i}{(i+1)(i+2)}
 &=\sum_{k=1}^{m}\frac1k\sum_{i=k}^{m}\left(\frac1{i+1}-\frac1{i+2}\right)\\
 &=\sum_{k=1}^{m}\frac1{k(k+1)}-\frac{H_m}{m+2}\\
 &=1-\frac1{m+1}-\frac{H_m}{m+2}.
\end{align*}
This proves~\eqref{eq:finiteN2}, including $m=0$. For model $O$,
$\ell_n^O(2)=4H_{n-1}/(n-1)-8/n$. Summing it and using
\[
 \sum_{j=1}^n\frac{H_j}{j}=\tfrac12(H_n^2+H_n^{(2)})
\]
gives~\eqref{eq:finiteO2}.

The root-variance checks are
\begin{align}
 \Var S_n^N(0)&=14-\frac{10}n-4\bigl(H_{n-1}^{(2)}+H_n^{(2)}\bigr),\label{eq:rootvarN}\\
 \Var S_n^O(0)&=2H_n-4H_n^{(2)}+2.\label{eq:rootvarO}
\end{align}
Both follow by summing $p_{j1}(1-p_{j1})$ over $2\le j\le n$; the root attachments are independent. Formula~\eqref{eq:rootvarN} agrees with~\cite[Proposition 4.1]{Lyon2020}.

\subsection{Exact finite-size second-moment recurrence}\label{app:momentrec}
Fix a truncation level $b\ge1$ and collect the state into a column vector
\[
 Y_n=(L_n(0),\ldots,L_n(b),W_n(0),\ldots,W_n(b))^{\mathsf T}.
\]
The symbol $Y_n$ is used only for this computational state. Let
$m_n=\E Y_n$ and $B_n=\E[Y_nY_n^{\mathsf T}]$. For $1\le r\le b$, let $u_{n,r}$ have entry $1$ in the $L(r)$ coordinate, entry $n$ in the $W(r)$ coordinate, and zeros elsewhere. The conditional probability of the update $u_{n,r}$ is a linear form $c_{n,r}^{\mathsf T}Y_{n-1}$, where
\begin{align*}
 c_{n,r}^{\mathsf T}Y_{n-1}
 &=\frac{2W_{n-1}(r-1)}{n(n-1)} &&\text{in model }N,\\
 &=\frac{2L_{n-1}(r-1)}{n-1}-\frac{2W_{n-1}(r-1)}{n(n-1)} &&\text{in model }O,\\
 &=\frac{L_{n-1}(r-1)}{n-1} &&\text{in model }U.
\end{align*}
An insertion at a larger depth leaves this truncated state unchanged. Since the possible depth updates are mutually exclusive,
\begin{equation}\label{eq:compmean}
 m_n=m_{n-1}+\sum_{r=1}^bu_{n,r}c_{n,r}^{\mathsf T}m_{n-1},
\end{equation}
while the uncentered second moments satisfy
\begin{align}
 B_n=B_{n-1}+\sum_{r=1}^b\bigl(&B_{n-1}c_{n,r}u_{n,r}^{\mathsf T}
       +u_{n,r}c_{n,r}^{\mathsf T}B_{n-1}\notag\\
       &+u_{n,r}u_{n,r}^{\mathsf T}c_{n,r}^{\mathsf T}m_{n-1}\bigr).\label{eq:compsecond}
\end{align}
The initial state has $L_1(0)=W_1(0)=1$ and all other coordinates zero, so $B_1=m_1m_1^{\mathsf T}$. The covariance matrix is $B_n-m_nm_n^{\mathsf T}$. In particular, for $b\ge2$ this computes the finite-size variance of $S_n(1)=L_n(1)+L_n(2)$ without a Gaussian approximation.

\section{A direct second-moment check for the uniform benchmark}\label{app:uniform}
Here all variables have the uniform law, and the superscript $U$ is suppressed. Let $X_n(r)=L_n(r)-\E L_n(r)$ and let $\varepsilon_n(r)$ be the centered conditional insertion indicator. The recurrence is
\begin{equation}\label{eq:UXrec}
 X_n(r)=\sum_{j=2}^n\frac{X_{j-1}(r-1)}{j-1}
             +\sum_{j=2}^n\varepsilon_j(r).
\end{equation}
The sum of martingale differences has squared norm at most $\E L_n(r)=O(\tau_n^r)$. Induction in~\eqref{eq:UXrec} then gives $\norm{X_n(r)}=O(\tau_n^{r-1/2})$: the drift has norm at most a constant times $\sum_{j\le n}\tau_j^{r-3/2}/j=O(\tau_n^{r-1/2})$, and the noise norm is $O(\tau_n^{r/2})$.

The independent root increments are
\[
 \xi_j^U=\1_{\{\Gamma_j=1\}}-\frac1{j-1}.
\]
Define $\alpha_n^U(1)=\sum_{j=2}^n\xi_j^U$ and
$\alpha_n^U(r)=\sum_{j=2}^n\alpha_{j-1}^U(r-1)/(j-1)$.
Subtracting this recurrence from~\eqref{eq:UXrec} gives
\[
 X_n(r)-\alpha_n^U(r)
 =\sum_{j=2}^n\frac{X_{j-1}(r-1)-\alpha_{j-1}^U(r-1)}{j-1}
       +\sum_{j=2}^n\varepsilon_j(r).
\]
At $r=2$ the first term is zero and the second has norm $O(\tau_n)$. For $r\ge3$, induction bounds the first term by $O(\tau_n^{r-1})$ and the second by $O(\tau_n^{r/2})$. Consequently the total error is $O(\tau_n^{r-1})$.

Iteration and Lemma~\ref{lem:kernel}, now without the factors $2^{r-1}$, give
\begin{equation}\label{eq:Urep}
 X_n(r)=\frac1{(r-1)!}\sum_{j=2}^n\xi_j^U\log^{r-1}(n/j)+R_n^U(r),
\end{equation}
where $R_n^U(1)=0$ and $\norm{R_n^U(r)}=O(\tau_n^{r-1})$ for $r\ge2$. Since
\[
 \Var(\xi_j^U)=\frac1{j-1}-\frac1{(j-1)^2}=\frac1j+O(j^{-2}),
\]
integral comparison yields
\[
 \sum_{j=2}^n\Var(\xi_j^U)\log^q(n/j)
 =\frac{(\log n)^{q+1}}{q+1}+O(\tau_n^q).
\]
Computing the covariance of the independent-increment sums in~\eqref{eq:Urep} gives the leading constant in~\eqref{eq:ucov}. Cauchy--Schwarz bounds a cross term involving one nonzero remainder by $O(\tau_n^{r+s-3/2})$ and the product of two remainders by $O(\tau_n^{r+s-2})$. These orders are smaller than $(\log n)^{r+s-1}$, completing the claimed second-moment verification. This establishes the covariance comparison used in Section~\ref{sec:clock}.

\section*{Data and code availability}\label{sec:code}
The code, data and figures are available in the public repository.\par
\noindent\mbox{\small\url{https://github.com/mgalazka84/fixed-level-age-dependent-recursive-trees}}\par
\noindent At \href{https://github.com/mgalazka84/fixed-level-age-dependent-recursive-trees/tree/71cefa4832416cf61265ba392f7263ea9ca9084e}{commit \texttt{71cefa483241}}, run \href{https://github.com/mgalazka84/fixed-level-age-dependent-recursive-trees/blob/71cefa4832416cf61265ba392f7263ea9ca9084e/code/check_kernel.py}{\texttt{python3 code/check\_kernel.py}} to reproduce the three intervals in~\eqref{eq:enclosures}. The script \texttt{code/check\_small\_trees.py} checks the moment recurrences by exact enumeration through seven vertices. The recorded simulation and both figures come from \href{https://github.com/mgalazka84/fixed-level-age-dependent-recursive-trees/tree/1765e49929ef1e14664fc887e10be2a86503dd9d}{commit \texttt{1765e49929ef}}. The repository includes the simulation seeds and file checksums.

\section*{Author contributions}
Marek Ga\l\k{a}zka: Conceptualization, Methodology, Formal analysis, Investigation, Validation, Writing -- original draft, Writing -- review and editing, Project administration.
Jerzy Szyma\'nski: Methodology, Formal analysis, Validation, Writing -- review and editing.

\section*{Funding and competing interests}
This research did not receive any specific grant from funding agencies in the public, commercial or not-for-profit sectors. The authors declare no competing financial interests or personal relationships that could have influenced the work.

\begingroup
\renewcommand{\bibfont}{\small}
\setlength{\bibsep}{4pt}

\endgroup
\section*{AI assistance with code}\label{sec:ai}
We used ChatGPT (OpenAI) to help write the simulation code and the scripts that verify the finite-size recurrences and the covariance enclosures in~\eqref{eq:enclosures}.
\end{document}